\documentclass[A4, 11pt]{article}

\usepackage{amssymb}
\usepackage{amsmath,amssymb,epsf,wrapfig,multicol,epic,eepic,trig,mathrsfs,dsfont,color,graphicx}
 \usepackage[abs]{overpic}

\def\qed{\hfill$\square$}
\newtheorem{thm}{Theorem}[section]
\newtheorem{cor}[thm]{Corollary}
\newtheorem{lem}[thm]{Lemma}
\newtheorem{prop}[thm]{Proposition}

\newtheorem{definition}[thm]{Definition}
\newtheorem{rem}[thm]{Remark}
\newtheorem{example}[thm]{Example}

\newcommand{\R}{\mathbb R}
\newcommand{\Z}{\mathbb Z}
\newcommand{\VLD}{{\cal VLD}}
\newcommand{\TLD}{{\cal TLD}}
\newcommand{\ALD}{{\cal ALD}}
\newcommand{\OALD}{{\cal OALD}}
\newcommand{\VL}{{\cal VL}}
\newcommand{\TL}{{\cal TL}}
\newcommand{\AL}{{\cal AL}}
\newcommand{\OAL}{{\cal OAL}}

\def\rmoveio#1#2{
\setlength{\unitlength}{#1}
\begin{picture}(50,30)
\put(5,0){\line(0,1){30}}

{\allinethickness{.8pt}
\put(10,15){\vector(1,0){15}}
\put(25,15){\vector(-1,0){15}}}

\qbezier(30,0)(30,20)(45,20)
\qbezier(45,20)(50,20)(50,15)
\qbezier(50,15)(50,10)(45,10)
\qbezier(45,10)(40,10)(36,14)
\qbezier(33,17)(30,25)(30,30)

\ifnum#2=2
\put(3,28){\path(0,0)(2,2)(4,0)}
\put(28,28){\path(0,0)(2,2)(4,0)}
\put(38,15){\makebox{${\Huge c_{1}}$}}
\fi

\end{picture}
}

\def\rmoveiio#1#2{
\setlength{\unitlength}{#1}
\begin{picture}(60,40)
\put(5,0){\line(0,1){40}}
\put(15,0){\line(0,1){40}}

{\allinethickness{.8pt}
\put(20,20){\vector(1,0){15}}
\put(35,20){\vector(-1,0){15}}}

\qbezier(40,0)(42,1)(47,5)
\qbezier(50,8)(68,20)(50,32)
\qbezier(47,35)(42,39)(40,40)

\qbezier(60,0)(20,20)(60,40)

\ifnum#2=2
\put(2,37){\path(0,0)(3,3)(6,0)}
\put(12,37){\path(0,0)(3,3)(6,0)}
\put(40,37){\path(0,0)(0,3)(3,3)}
\put(60,37){\path(0,0)(0,3)(-3,3)}
\put(47,25){\makebox{${\Huge c_{1}}$}}
\put(47,13){\makebox{${\Huge c_{2}}$}}
\fi

\ifnum#2=3
\put(2,2){\path(0,0)(3,-3)(6,0)}
\put(12,37){\path(0,0)(3,3)(6,0)}
\put(40,3){\path(0,0)(0,-3)(3,-3)}
\put(60,37){\path(0,0)(0,3)(-3,3)}
\put(47,25){\makebox{${\Huge c_{1}}$}}
\put(47,13){\makebox{${\Huge c_{2}}$}}
\fi

\end{picture}
}
\def\rmoveiiio#1#2{
\setlength{\unitlength}{#1}
\begin{picture}(75,30)
\put(0,0){\line(1,1){15}}
\qbezier(15,15)(20,20)(20,30)

\put(10,0){\line(-1,1){4}}
\qbezier(4,6)(-5,15)(5,25)
\put(5,25){\line(1,1){5}}

\qbezier(20,0)(20,10)(16,14)
\put(14,16){\line(-1,1){8}}
\put(4,26){\line(-1,1){4}}

{\allinethickness{.8pt}
\put(30,15){\vector(1,0){10}}
\put(40,15){\vector(-1,0){10}}}

\qbezier(50,0)(50,10)(55,15)
\put(55,15){\line(1,1){15}}

\put(60,0){\line(1,1){5}}
\qbezier(65,5)(75,15)(66,24)
\put(64,26){\line(-1,1){4}}

\put(70,0){\line(-1,1){4}}
\put(64,6){\line(-1,1){8}}
\qbezier(54,16)(50,20)(50,30)

\ifnum#2=2
\put(0,27){\path(0,0)(0,3)(3,3)}
\put(7,30){\path(0,0)(3,0)(3,-3)}
\put(17,27){\path(0,0)(3,3)(6,0)}

\put(10,23){\makebox{${\Huge c_{1}}$}}
\put(10,3){\makebox{${\Huge c_{2}}$}}
\put(20,13){\makebox{${\Huge c_{3}}$}}

\put(47,27){\path(0,0)(3,3)(6,0)}
\put(60,27){\path(0,0)(0,3)(3,3)}
\put(67,30){\path(0,0)(3,0)(3,-3)}

\put(70,23){\makebox{${\Huge c'_{2}}$}}
\put(70,3){\makebox{${\Huge c'_{1}}$}}
\put(60,13){\makebox{${\Huge c'_{3}}$}}
\fi

\end{picture}
}
\def\rmovevio#1#2{
\setlength{\unitlength}{#1}
\begin{picture}(50,30)
\put(5,0){\line(0,1){30}}

{\allinethickness{.8pt}
\put(10,15){\vector(1,0){15}}
\put(25,15){\vector(-1,0){15}}}

\qbezier(30,0)(30,20)(45,20)
\qbezier(45,20)(50,20)(50,15)
\qbezier(50,15)(50,10)(45,10)
\qbezier(45,10)(30,10)(30,30)

\put(34,15){\circle{5}}

\ifnum#2=2
\put(3,28){\path(0,0)(2,2)(4,0)}
\put(28,28){\path(0,0)(2,2)(4,0)}
\put(38,15){\makebox{${\Huge c_{1}}$}}
\fi

\end{picture}
}

\def\rmoveviio#1#2{
\setlength{\unitlength}{#1}
\begin{picture}(60,40)
\put(5,0){\line(0,1){40}}
\put(15,0){\line(0,1){40}}

{\allinethickness{.8pt}
\put(20,20){\vector(1,0){15}}
\put(35,20){\vector(-1,0){15}}}

\qbezier(40,0)(80,20)(40,40)
\qbezier(60,0)(20,20)(60,40)
\put(50,6){\circle{5}}
\put(50,34){\circle{5}}

\ifnum#2=2
\put(2,37){\path(0,0)(3,3)(6,0)}
\put(12,37){\path(0,0)(3,3)(6,0)}
\put(40,37){\path(0,0)(0,3)(3,3)}
\put(60,37){\path(0,0)(0,3)(-3,3)}
\put(47,25){\makebox{${\Huge c_{1}}$}}
\put(47,13){\makebox{${\Huge c_{2}}$}}
\fi

\ifnum#2=3
\put(2,2){\path(0,0)(3,-3)(6,0)}
\put(12,37){\path(0,0)(3,3)(6,0)}
\put(40,3){\path(0,0)(0,-3)(3,-3)}
\put(60,37){\path(0,0)(0,3)(-3,3)}
\put(47,25){\makebox{${\Huge c_{1}}$}}
\put(47,13){\makebox{${\Huge c_{2}}$}}
\fi

\end{picture}
}

\def\rmoveviiio#1#2{
\setlength{\unitlength}{#1}
\begin{picture}(70,30)
\put(0,0){\line(1,1){15}}
\qbezier(15,15)(20,20)(20,30)

\put(10,0){\line(-1,1){5}}
\qbezier(5,5)(-5,15)(5,25)
\put(5,25){\line(1,1){5}}

\qbezier(20,0)(20,10)(15,15)
\put(15,15){\line(-1,1){15}}

\put(5,5){\circle{5}}
\put(15,15){\circle{5}}
\put(5,25){\circle{5}}

{\allinethickness{.8pt}
\put(30,15){\vector(1,0){10}}
\put(40,15){\vector(-1,0){10}}}

\qbezier(50,0)(50,10)(55,15)
\put(55,15){\line(1,1){15}}

\put(60,0){\line(1,1){5}}
\qbezier(65,5)(75,15)(65,25)
\put(65,25){\line(-1,1){5}}

\put(70,0){\line(-1,1){15}}
\qbezier(55,15)(50,20)(50,30)

\put(65,5){\circle{5}}
\put(55,15){\circle{5}}
\put(65,25){\circle{5}}

\ifnum#2=2
\put(0,27){\path(0,0)(0,3)(3,3)}
\put(7,30){\path(0,0)(3,0)(3,-3)}
\put(17,27){\path(0,0)(3,3)(6,0)}

\put(20,13){\makebox{${\Huge c_{1}}$}}

\put(47,27){\path(0,0)(3,3)(6,0)}
\put(60,27){\path(0,0)(0,3)(3,3)}
\put(67,30){\path(0,0)(3,0)(3,-3)}

\put(60,13){\makebox{${\Huge c'_{1}}$}}
\fi

\end{picture}
}

\def\rmovevivo#1#2{
\setlength{\unitlength}{#1}
\begin{picture}(70,30)
\put(0,0){\line(1,1){15}}
\qbezier(15,15)(20,20)(20,30)

\put(10,0){\line(-1,1){5}}
\qbezier(5,5)(-5,15)(5,25)
\put(5,25){\line(1,1){5}}

\qbezier(20,0)(20,10)(16,14)
\put(14,16){\line(-1,1){14}}

\put(5,5){\circle{5}}
\put(5,25){\circle{5}}

{\allinethickness{.8pt}
\put(30,15){\vector(1,0){10}}
\put(40,15){\vector(-1,0){10}}}

\qbezier(50,0)(50,10)(55,15)
\put(55,15){\line(1,1){15}}

\put(60,0){\line(1,1){5}}
\qbezier(65,5)(75,15)(65,25)
\put(65,25){\line(-1,1){5}}

\put(70,0){\line(-1,1){14}}
\qbezier(54,16)(50,20)(50,30)

\put(65,5){\circle{5}}
\put(65,25){\circle{5}}

\ifnum#2=2
\put(0,27){\path(0,0)(0,3)(3,3)}
\put(7,30){\path(0,0)(3,0)(3,-3)}
\put(17,27){\path(0,0)(3,3)(6,0)}

\put(20,13){\makebox{${\Huge c_{1}}$}}

\put(47,27){\path(0,0)(3,3)(6,0)}
\put(60,27){\path(0,0)(0,3)(3,3)}
\put(67,30){\path(0,0)(3,0)(3,-3)}

\put(60,13){\makebox{${\Huge c'_{1}}$}}
\fi

\end{picture}
}
\def\rmoveti#1{
\setlength{\unitlength}{#1}
\begin{picture}(60,20)
\put(0,10){\line(1,0){20}}
\put(15,0){\line(0,1){20}}
\put(15,10){\circle{5}}
{\linethickness{2pt}
\put(7,7){\line(0,1){6}}} 
{\allinethickness{.8pt}
\put(25,10){\vector(1,0){10}}
\put(35,10){\vector(-1,0){10}}}

\put(40,10){\line(1,0){20}}
\put(45,0){\line(0,1){20}}
\put(45,10){\circle{5}}
{\linethickness{2pt}
\put(53,7){\line(0,1){6}}}
\end{picture}
}

\def\rmovetiio#1{
\setlength{\unitlength}{#1}
\begin{picture}(40,20)
\put(5,0){\line(0,1){20}}
{\linethickness{2pt}
\put(2.,7){\line(1,0){6}} 
\put(2.,13){\line(1,0){6}} }

{\allinethickness{.8pt}
\put(15,10){\vector(1,0){10}}
\put(25,10){\vector(-1,0){10}}}
\put(35,0){\line(0,1){20}}
\end{picture}
}

\def\rmovetiiio#1#2{
\setlength{\unitlength}{#1}
\begin{picture}(100,30)

\put(0,0){\line(1,1){20}}
\put(20,0){\line(-1,1){9}}
\put(9,11){\line(-1,1){9}}

{\allinethickness{2pt}
\put(7,3){\line(-1,1){4}} 
\put(13,3){\line(1,1){4}} 
\put(3,13){\line(1,1){4}} 
\put(17,13){\line(-1,1){4}} }

{\allinethickness{.8pt}
\put(25,10){\vector(1,0){10}}
\put(35,10){\vector(-1,0){10}}}

\qbezier(40,2)(48,20)(54,11)
\qbezier(56, 9)(62,0)(70,18)

\qbezier(40,18)(48,0)(55,10)
\qbezier(55,10)(62,20)(70,2)

\put(45,10){\circle{5}}
\put(65,10){\circle{5}}

\put(73,10){or}

\qbezier(80,-3)(98,5)(89,9)
\qbezier(87,11)(78,15)(96,23)

\qbezier(96,-3)(78,5)(88,10)
\qbezier(88,10)(98,15)(80,23)

\put(88,1){\circle{5}}
\put(88,19){\circle{5}}

\ifnum#2=2
\put(0,17){\path(0,0)(0,3)(3,3)}
\put(17,20){\path(0,0)(3,0)(3,-3)}

\put(15,8){\makebox{${\Huge c_{1}}$}}

\put(40,12){\path(0,0)(0,3)(3,3)}
\put(67,15){\path(0,0)(3,0)(3,-3)}

\put(53,3){\makebox{${\Huge c'_{1}}$}}
\fi

\end{picture}
}

\def\fgknotgpp#1{
\setlength{\unitlength}{#1}
\begin{picture}(50,50)
\put(25,45){\vector(0,-1){40}}
\put(5,25){\line(1,0){17}}
\put(28,25){\vector(1,0){17}}
\put(0,30){\makebox{$x_j$}}
\put(0,15){\makebox{$y_j$}}
\put(40,30){\makebox{$x_{j+1}$}}
\put(40,15){\makebox{$y_{j+1}$}}
\put(15,45){\makebox{$y_i$}}
\put(30,45){\makebox{$x_i$}}
\put(15,0){\makebox{$y_{i+1}$}}
\put(30,0){\makebox{$x_{i+1}$}}
\end{picture}
}
\def\fgknotgpm#1{
\setlength{\unitlength}{#1}
\begin{picture}(50,50)
\put(25,45){\vector(0,-1){40}}
\put(22,25){\vector(-1,0){17}}
\put(45,25){\line(-1,0){17}}
\put(0,30){\makebox{$y_{j+1}$}}
\put(0,15){\makebox{$x_{j+1}$}}
\put(40,30){\makebox{$y_{j}$}}
\put(40,15){\makebox{$x_{j}$}}
\put(15,45){\makebox{$y_i$}}
\put(30,45){\makebox{$x_i$}}
\put(15,0){\makebox{$y_{i+1}$}}
\put(30,0){\makebox{$x_{i+1}$}}
\end{picture}
}

\def\fgknotgpv#1{
\setlength{\unitlength}{#1}
\begin{picture}(50,50)
\put(25,45){\vector(0,-1){40}}
\put(5,25){\vector(1,0){40}}
\put(25,25){\circle{10}}
\put(0,30){\makebox{$x_j$}}
\put(0,15){\makebox{$y_j$}}
\put(40,30){\makebox{$x_j$}}
\put(40,15){\makebox{$y_j$}}
\put(15,45){\makebox{$y_i$}}
\put(30,45){\makebox{$x_i$}}
\put(15,0){\makebox{$y_i$}}
\put(30,0){\makebox{$x_i$}}
\end{picture}
}

\def\fgknotgpt#1{
\setlength{\unitlength}{#1}
\begin{picture}(20,50)
\put(10,45){\vector(0,-1){40}}
{\linethickness{2pt}
\put(7,25){\line(1,0){6}}}
\put(0,45){\makebox{$y_i$}}
\put(15,45){\makebox{$x_i$}}
\put(0,0){\makebox{$y_{i+1}$}}
\put(15,0){\makebox{$x_{i+1}$}}
\end{picture}
}

\begin{document}
\title
{Double coverings of twisted links\thanks{Dedicated to the honor of our late colleague Slavik Jablan.}} 
\author{
Naoko Kamada
\thanks{This work was supported by JSPS KAKENHI Grant Numbers 15K04879.}
\\ Graduate School of Natural Sciences,  Nagoya City University\\ 
Mizuho-ku, Nagoya, Aichi 467-8501 Japan\\
\and 
Seiichi Kamada
\thanks{This work was supported by JSPS KAKENHI Grant Numbers 26287013.}
\\ Department of Mathematics, Osaka City University, \\
Suimiyoshi,  Osaka 558-8585, Japan
}

\date{}

\maketitle

\begin{abstract} 
Twisted links are a generalization of virtual links.  As virtual links correspond to abstract links on orientable surfaces, twisted links correspond to abstract links on (possibly non-orientable)  surfaces.  In this paper, we introduce the notion of the double covering of a twisted link.   It is defined by considering the orientation double covering of an abstract link or alternatively by constructing a diagram called a double covering diagram.  We also discuss links in thickened surfaces, their diagrams and their stable equivalence classes.   Bourgoin's twisted knot group is understood as the virtual knot group of the double covering.  
\end{abstract}



\section{Introduction}\label{sect:intro}

L. Kauffman \cite{rkauD} introduced the notion of virtual knots and links, which is a generalization of classical knots and links. 
Virtual links are in one-to-one correspondence to 
abstract links on oriented surfaces \cite{rkk}, and to 
stable equivalence classes of  links in oriented thickened 
surfaces over oriented surfaces \cite{rCKS, rFRS}. 
M. Bourgoin  \cite{rBourg} generalized them to twisted links.  
Twisted links are in one-to-one correspondence to abstract links on (possibly non-orientable) surfaces, and  to stable equivalence classes of  links in oriented thickened surfaces over  (possibly non-orientable)  surfaces \cite{rBourg}. 
 
 In this paper, we introduce the notion of the double covering of a twisted link (Section~\ref{sect:doubleA}). 
 It is defined by considering the orientation double covering of an abstract link.  
 We  give an  explicit method of constructing a virtual link diagram $\tilde D$ from a twisted link diagram $D$, called a double covering diagram, such that the virtual link $\tilde K=[ \tilde D]$ is the double covering of the twisted link $K=[D]$ (Section~\ref{sect:doubleB}). 
 
As an application, one can obtain an invariant of twisted links from an invariant of virtual links by considering  the double coverings.  For example, Bourgoin's twisted link group \cite{rBourg} is understood as the virtual link group of the double covering. This gives a geometric interpretation of the twisted link group, which  answers  
the problem proposed by Bourgoin in \cite{rBourg}.  Similarly the twisted link  quandle of a twisted link defined by the first author in \cite{rkamF} is the (upper) virtual link quandle of the double covering.  These are discussed in Section~\ref{sect:groups}.  

We also discuss links in oriented thickened surfaces, their diagrams and stable equivalence classes in Sections~\ref{sect:thickenedA} and~\ref{sect:thickenedB}. 
There is a bijection between the set of stable equivalence classes of links in oriented thickened surfaces over closed surfaces (or equivalently the set of stable equivalence classes of generalized link diagrams over closed surfaces) and the set of abstract links.

\begin{figure}[ht]
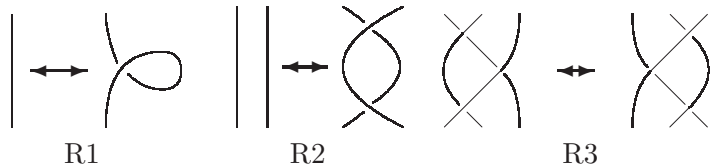

\begin{center}
\begin{tabular}{ccc}
\rmoveio{.5mm}{1}&\rmoveiio{.4mm}{1}&\rmoveiiio{.5mm}{1}\\
R1 & R2 & R3 \\
\end{tabular}
\caption{Reidemeister moves (R-moves)}\label{fig:movesR}
\end{center}
\end{figure}

\begin{figure}[ht]
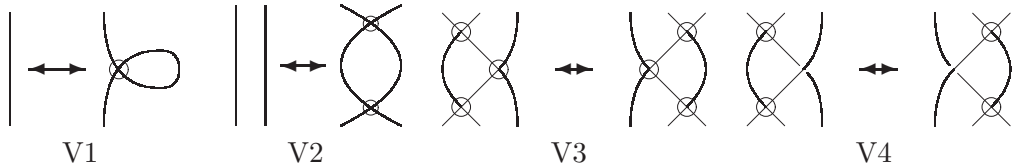

\begin{center}
\begin{tabular}{cccc}
\rmovevio{.5mm}{1}&\rmoveviio{.4mm}{1}&\rmoveviiio{.5mm}{1}&\rmovevivo{.5mm}{1}\\
V1 &  V2 &  V3 & V4 \\
\end{tabular}
\caption{Virtual Reidemeister moves (V-moves)}\label{fig:movesV}
\end{center}
\end{figure}

A {\it virtual link diagram\/} is a link diagram in $\R^2$ which may have virtual crossings;   
a  {\it virtual crossing\/} is an encircled double point without over-under information. 
A {\it virtual link\/} is the equivalence class of a virtual link diagram where the equivalence relation is generated by isotopies of $\R^2$, Reidemeister moves (Fig.~\ref{fig:movesR}) and virtual Reidemeister moves (Fig.~\ref{fig:movesV}).  

\begin{figure}[ht]
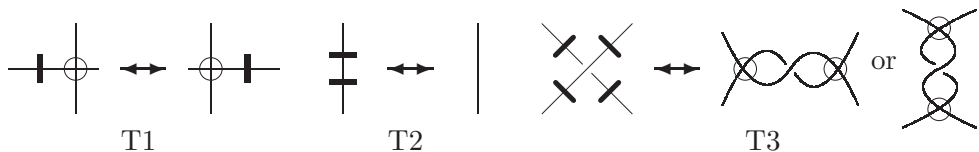

\begin{center}
\begin{tabular}{ccc}
\rmoveti{.6mm}&\rmovetiio{.6mm}&\rmovetiiio{.6mm}{1}\\
T1 & T2 & T3 \\
\end{tabular}
\caption{Twisted Reidemeister moves (T-moves)}\label{fig:movesT}
\end{center}
\end{figure}

A {\it twisted link diagram\/} is a link diagram which may have virtual crossings or 
bars; a {\it bar} is a small edge intersecting transversely with an arc of the diagram.  
A {\it twisted link\/} is the equivalence class of a twisted link diagram where the equivalence relation is generated by isotopies of $\R^2$, Reidemeister moves, virtual Reidemeister moves and twisted Reidemeister moves 
(Fig.~\ref{fig:movesT}).  Reidemeister moves, virtual Reidemeister moves and twisted Reidemeister moves are called R-moves, V-moves and T-moves for short. 
A virtual link (or a twisted link) with one component is called a {\it virtual knot} 
(or a {\it twisted knot}).

\vspace{0.5cm}

Throughout this paper, we assume that virtual links and twisted links are 
oriented, although orientations are not specified in some figures.  

\begin{figure}[ht]
\begin{center}
\includegraphics[width=11cm]{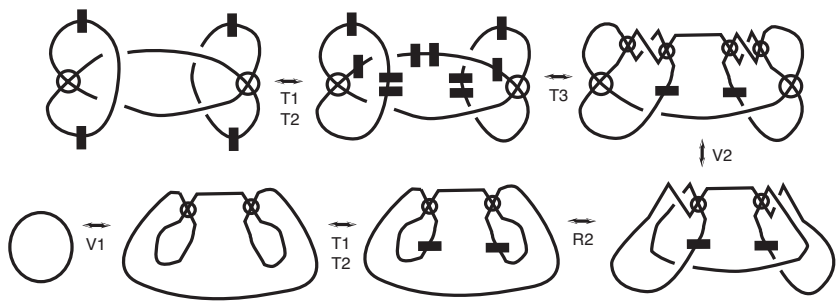}
\caption{A twisted knot diagram equivalent to the trivial diagram}\label{fig:figmoveA} 
\end{center}
\end{figure}

Here are some examples. 
In Fig.~\ref{fig:figmoveA} the top left is a diagram obtained by adding bars to a virtual knot diagram called Kishino's virtual knot diagram.  It is equivalent to the trivial knot diagram as a twisted knot as shown in the figure. 

\begin{figure}[ht]
\begin{center}
\includegraphics[width=11cm]{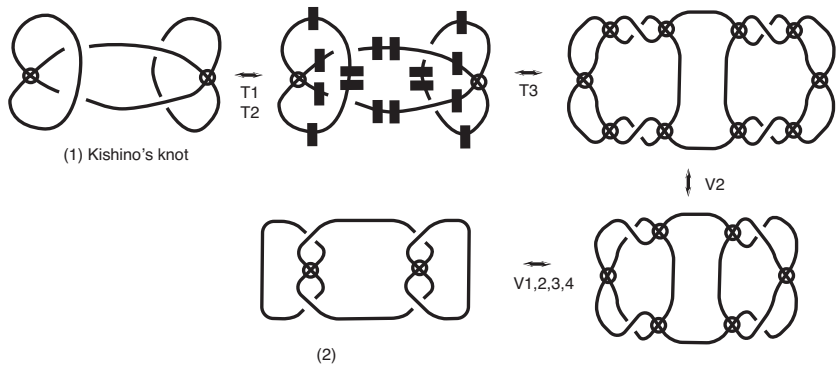}
\caption{Twisted knot diagrams equivalent as twisted knots}
\label{fig:figmoveB} 
\end{center}
\end{figure}

In Fig.~\ref{fig:figmoveB} the top left is Kishino's virtual knot diagram.  It is equivalent as a twisted knot to the diagram on the bottom left.    

\vspace{0.2cm}

The idea of the double covering of a twisted link was first announced by the authors at the conference \lq\lq Knots in Washington XXX'' held at George Washington University in May, 2010.  In this paper, we give the details. 

\section{The double covering of a twisted link}\label{sect:doubleA}

We introduce the notion of the double covering of a twisted link.  
The definition is given in terms of abstract links, and an explicit method of constructing a diagram representing the double covering will be given in the next section.  
 
Abstract links in the sense of \cite{rkk} are considered in the case where base surfaces are oriented, and they are generalized by Bourgoin  \cite{rBourg} into the case where base surfaces are not necessarily orientable.  First we recall them.  Our definition in this paper is a slightly modified version of the original definitions in  \cite{rkk} and \cite{rBourg}.   

An {\it abstract link diagram} is a link diagram $D$ on a compact surface 
$\Sigma$ 
such that the underlying $4$-valent graph $|D|$ is a deformation retract of $\Sigma$ and 
for each crossing point of the diagram a local orientation of the surface $\Sigma$ is assigned.  (A {\it local orientation} of $\Sigma$ at a point $x \in \Sigma$ means a generator of $H_2(N(x), \partial N(x); \Z) \cong H_2(N(x), N(x) \setminus \{x\}; \Z) \cong 
H_2(\Sigma, \Sigma \setminus \{x\}; \Z)$ where $N(x)$ is a regular neighborhood of $x$. 
An over-under information on a vertex $x$ of $|D|$ is a choice of a pair of edges around $v$ in a diagonal position as `over' and the other pair as `under'. Such a choice is expressed by removing a small portion from the `under' edges to make a crossing of $D$.) 
We call $\Sigma$ the {\it base surface}.  The collection $O$ of the local orientations assigned to the crossings of $D$ is called the {\it local orientation system} of the abstract link diagram. 

An abstract link diagram $D$ on $\Sigma$ with local orientation system $O$ is  denoted by the triple $(\Sigma, D, O)$.   

Two abstract link diagrams $(\Sigma, D, O)$ and $(\Sigma', D', O')$ are identified or considered to be the same if there is a homeomorphism from $\Sigma$ to $\Sigma'$ carrying $|D|$ to $|D'|$ such that for each crossing point of $D$, the 
local orientation and the over-under information are preserved.  

Let $(\Sigma, D, O)$ be an abstract link diagram and let $x$ be a crossing point of $D$. 
An abstract link diagram $(\Sigma', D', O')$ is said to be obtained from $(\Sigma, D, O)$ by an  {\it inversion} at $x$ if it is obtained by reversing the local orientation at $x$ and reversing the over-under information of the crossing $x$. 

An abstract link diagram $(\Sigma', D', O')$ is said to be obtained from $(\Sigma, D, O)$ by an {\it abstract R-move} if there exist two embeddings  
$f: \Sigma \to F$ and $f': \Sigma' \to F$ to a surface $F$ and there exists an oriented $2$-disk $N$ in $F$ such that 
the diagram $f(D')$ is obtained from $f(D)$ by a Reidemeister move on $F$  in 
$N$ and the local orientations of the crossings of $f(D)$ and $f'(D)$ in $N$ are compatible with the orientation of $N$.  (For a point $x \in N$, a local orientation of $N$ at $x$ is {\it compatible with} the orientation of $N$ if they coincide under the isomorphism 
$H_2(N, \partial N; \Z) \cong H_2(N, N \setminus \{x\}; \Z)$ induced from the inclusion map.) 

\begin{definition}{\rm 
Two abstract link diagrams $A= (\Sigma, D, O)$ and $A'= (\Sigma', D', O')$ are 
 {\it equivalent} if there exists a finite sequence of abstract link diagrams 
$A =  A_0,  A_1, \dots,  A_n = A'$ such that for each $i \in \{1, \dots, n\}$, $A_i$ is obtained from $A_{i-1}$ by 
an inversion or an abstract R-move.   
An {\it abstract link} is the equivalence class of an abstract link diagram.  
}\end{definition}

\begin{definition}{\rm 
An {\it ordinary abstract link diagram} or {\it abstract link diagram in the ordinary sense} is an abstract link diagram $A=(\Sigma, D, O)$ such that 
$\Sigma$ is oriented and the local orientation at every crossing point is compatible with the orientation of $\Sigma$.  In this case, we may denote $A=(\Sigma, D, O)$ by the pair $A=(\Sigma, D)$, since there is no need to specify $O$. 
An {\it ordinary abstract link} is an abstract link represented by an ordinary abstract link diagram.  
}\end{definition} 

\begin{figure}[ht]
\begin{center}
\includegraphics[width=10cm]{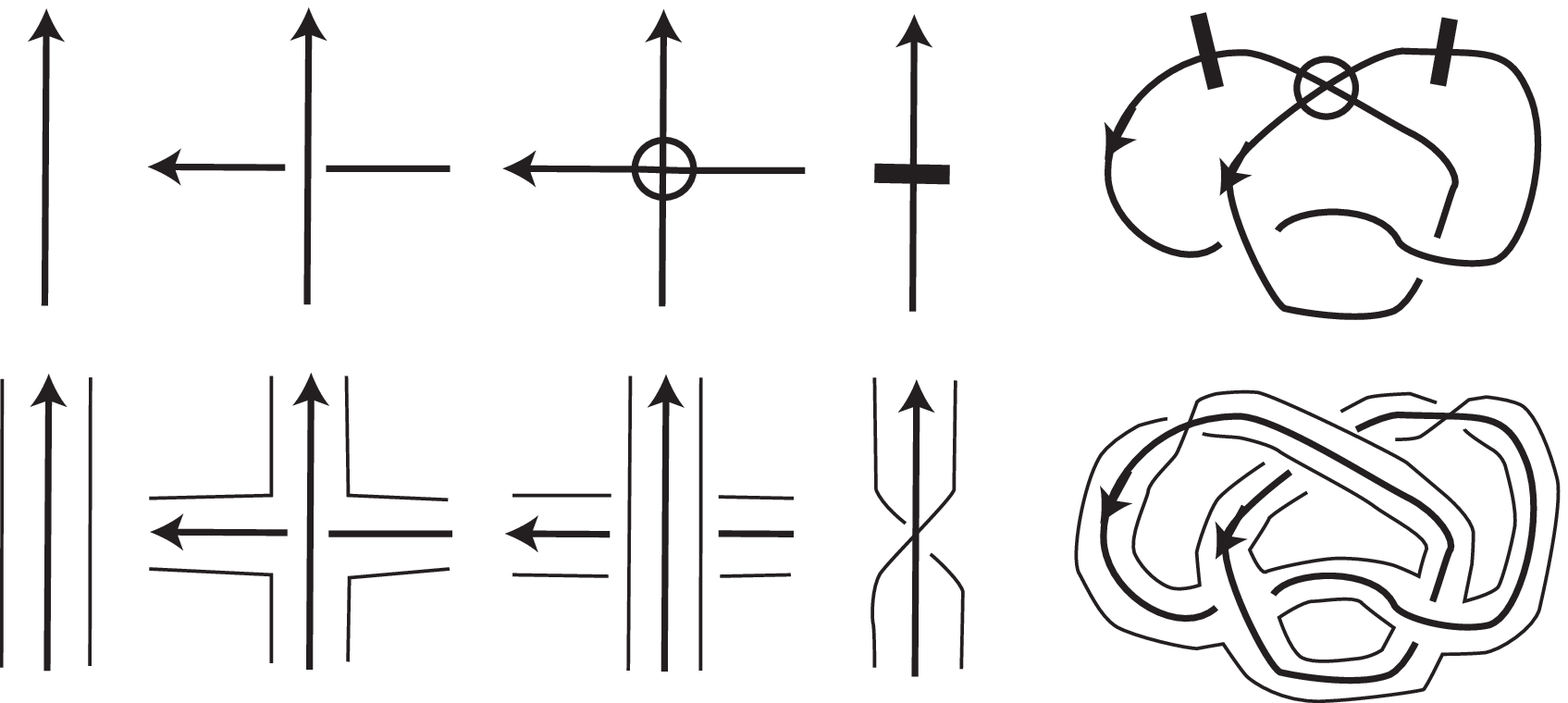}
\caption{A twisted link diagram and its associated abstract link diagram}
\label{fig:figabst} 
\end{center}
\end{figure}

\begin{figure}[ht]
\begin{center}
\includegraphics[width=7cm]{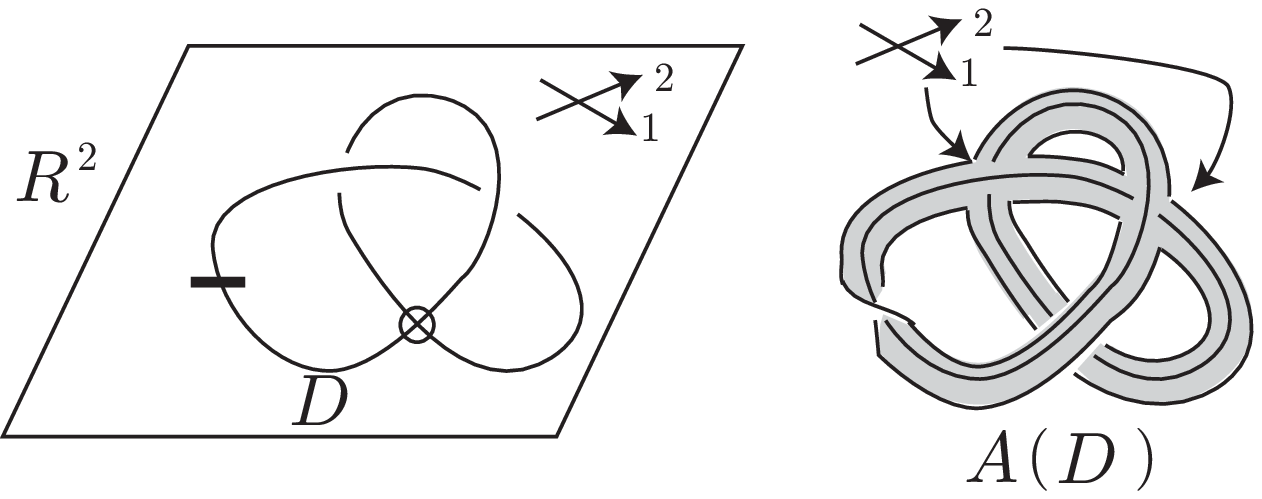}
\caption{A twisted link diagram and its associated abstract link diagram}
\label{fig:figdiagd} 
\end{center}
\end{figure}

The authors \cite{rkk} introduced a method of constructing an ordinary abstract link diagram from a virtual link diagram, and Bourgoin \cite{rBourg} generalized it to a method of constructing an abstract link diagram from a twisted link diagram.  
The method is illustrated in Fig.~\ref{fig:figabst}.  
For a twisted link diagram $D$, we denote by $A(D)= (\Sigma, A(D), O)$ the 
abstract link diagram obtained by this method and call it 
the {\it abstract link diagram associated with $D$}. 
 (We abuse the symbol $A(D)$ for the diagram itself and for the abstract link.) 
Here we assume that the local orientation 
of $\Sigma$ at each crossing point of $A(D)$ is the orientation inherited from the standard orientation of $\R^2$.  
See Fig.~\ref{fig:figdiagd}.  
For the precise definition and details, refer to 
\cite{rkk} for a virtual link diagram and \cite{rBourg} for a twisted link diagram. 

\vspace{0.5cm} 

\vspace{0.5cm} 

We denote by $\ALD$, $\OALD$, $\VLD$ and $\TLD$ the set of 
abstract link diagrams, the set of ordinary abstract link diagrams, the set of virtual link diagrams and the set of twisted link diagrams, and denote by 
$\AL$, $\OAL$, $\VL$ and $\TL$ the set of 
abstract links, the set of ordinary abstract links, the set of virtual links and the set of twisted links, respectively. 

\begin{thm}[Kamada-Kamada \cite{rkk} and Bourgoin \cite{rBourg}]\label{thm:Abst} 
The map 
$\varphi_{\rm ord}:  \VLD \to \OALD, D \mapsto A(D)$ 
induces a bijection $\Phi_{\rm ord}:  \VL \to \OAL, [D] \mapsto [A(D)]$.  
The map 
$\varphi : \TLD \to \ALD, D \mapsto A(D)$ 
induces a bijection  
$\Phi: \TL \to \AL, [D] \mapsto [A(D)]$.  
\end{thm}

A double covering of a surface $F$ is called the {\it orientation double covering} of $F$ if   the restriction over each connected orientable  component of $F$ is the trivial double covering and 
the restriction over each connected non-orientable component is a double covering  
whose total space is orientable.  

We introduce the double covering of an abstract link diagram.  

\begin{definition}{\rm 
An abstract link diagram $A' = (\Sigma', D', O')$ is the 
{\it orientation double covering} of $A = (\Sigma, D, O)$ if there is a map 
$\pi :  \Sigma' \to \Sigma$ satisfying the following. 
\begin{itemize}
\vspace{-0.2cm}
\item[(1)] $\pi :  \Sigma' \to \Sigma$ is the orientation double covering of $\Sigma$.  
\vspace{-0.2cm}
\item[(2)] The restriction to the underlying $4$-valent graph $|D'|$ is a double covering 
$|D'| \to |D|$.
\vspace{-0.2cm}
\item[(3)] For each crossing $x' $ of $D'$, the over-under information at $x'$ of $D'$ is mapped to the over-under information at $x = \pi(x')$ of $D$. 
\vspace{-0.2cm}
\item[(4)] For each crossing $x' $ of $D'$, the local orientation at $x'$ of $O'$ is mapped to the local orientation at $x = \pi(x')$ of $O$. 
\end{itemize}
}\end{definition}

Now we define the double covering of a twisted link.  

Let $K$ be a twisted link and let $D$ be a twisted link diagram representing $K$.  
Let $A(D)= (\Sigma, A(D), O)$ be the abstract link diagram associated with $D$, and let 
$A'(D)= (\Sigma', A'(D), O')$ be the orientation double covering of $A(D)$.  Note that $\Sigma'$ is an orientable surface.  Give an orientation to $\Sigma'$ and denote by $\Sigma''$ the oriented surface.  Apply an inversion to  $A'(D)= (\Sigma', A'(D), O')$ at every crossing point $x'$ of $A'(D)$ such that the local orientation at $x'$ is not compatible with the orientation of $\Sigma''$.  We obtain an ordinary abstract link diagram $A''(D) = (\Sigma'', A''(D))$.  By Theorem~\ref{thm:Abst}, there is a virtual link, say  $\tilde K$, such that 
$\Phi_{\rm ord}(\tilde K) = [A''(D)]$.  

\begin{definition}{\rm 
In the above situation, the virtual link  $\tilde K$ is called the {\it orientation double covering} or simply the   
{\it double covering} of $K$.  
}\end{definition} 

\begin{figure}[ht]
\begin{center}
\includegraphics[width=10cm]{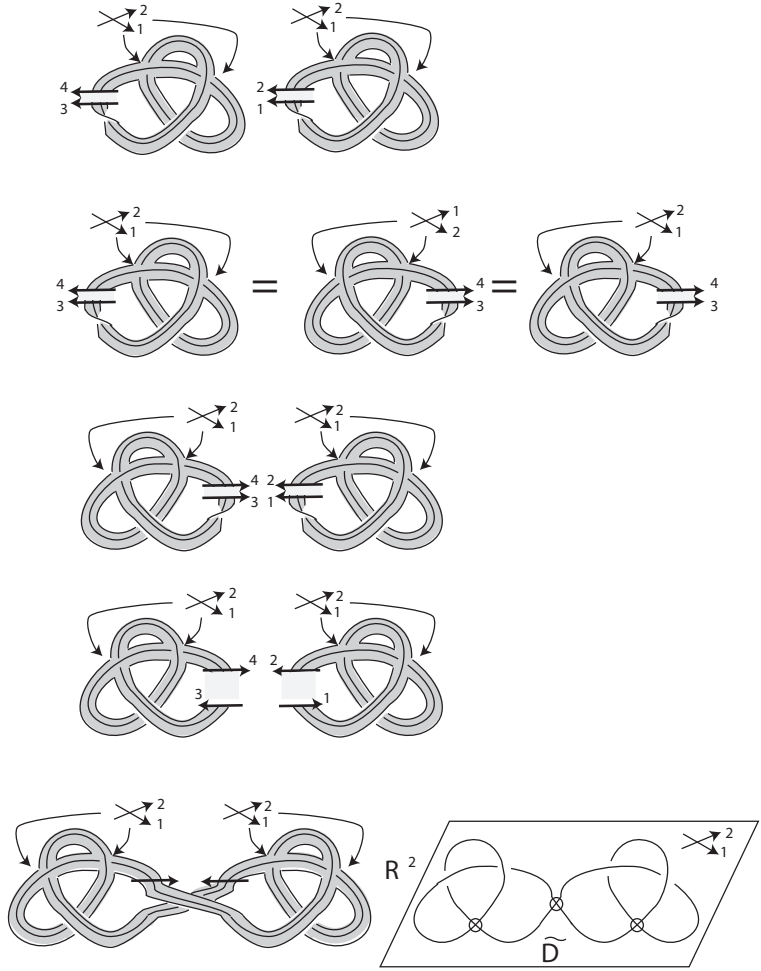}
\caption{Orientation double covering of $A(D)$}
\label{fig:figdiagfa} 
\end{center}
\end{figure}

\begin{example}\label{example:A}{\rm  
Let $K$ be the twisted link represented by 
the twisted link diagram $D$ depicted in Fig.~\ref{fig:figdiagd}.  
Let $A(D)$ be the abstract link diagram depicted in Fig.~\ref{fig:figdiagd}, which is 
the abstract link diagram associated with $D$.   
We consider the orientation double covering of $A(D)$, modulo inversions.  
First we prepare two copies of $A(D)$ and cut them as in the top row of Fig.~\ref{fig:figdiagfa}.  Identify the arrows labeled $1$ and $4$ and identify the arrows labeled $2$ and $3$, and we obtain the orientation double covering $A'(D)$ of $A(D)$.  

The left side of the top row of the figure can be changed as shown in the second row:   
The first equality is obtained by considering the reflection $(x,y,z) \mapsto (-x,y,z)$ when we consider them in the $3$-space.  The second equality is obtained by inversions.  
Thus, when we  identify the arrows labeled $1$ and $4$ and identify the arrows labeled $2$ and $3$, 
we see that the union of the two abstract link diagrams in the third row of Fig.~\ref{fig:figdiagfa} is the orientation double covering of $A(D)$ modulo inversions.  

As in the fourth row, we remove the twists of the surfaces and finally obtain an abstract link diagram as in the bottom left of Fig.~\ref{fig:figdiagfa}.  This is the orientation double covering of $A(D)$ modulo   inversions, say $A''(D)$.  
Note that this abstract link diagram is the ordinary abstract link diagram $A(\tilde D)$ associated with 
the virtual link diagram $\tilde D$ depicted in the figure.  
Thus the double covering of $K$ is the virtual link represented by 
$\tilde D$.  (The diagram $\tilde D$ is discussed in the next section.)
}\end{example}

\vspace{0.5cm} 

\begin{thm}\label{thm:doubleA}
The double covering of a twisted link is well defined.  
\end{thm} 

Proof. 
Suppose that $D_1$ and $D_2$ represent the same twisted link $K$.  
By Theorem~\ref{thm:Abst}, 
the abstract link diagrams $A(D_1)$ and $A(D_2)$ 
associated with $D_1$ and $D_2$ are equivalent as abstract links.  
Their orientation double coverings $A'(D_1)$ and $A'(D_2)$ are equivalent as abstract links. 
(If $A(D_2)$ is obtained from $A(D_1)$ by an inversion on a crossing point $x$, then  $A'(D_2)$ is obtained from $A'(D_1)$ by inversions at the two crossings in the preimage of $x$. 
If $A(D_2)$ is obtained from $A(D_1)$ by an abstract R-move, then $A'(D_2)$ is obtained from $A'(D_1)$ by a pair of abstract R-moves that cover the abstract R-move.) 
Thus ordinary abstract link diagrams $A''(D_1)$ and $A''(D_2)$ are equivalent. 
By Theorem~\ref{thm:Abst},  we obtain a unique virtual link as the double covering of $K$.  
\qed

\section{Construction of a diagram of the double covering}\label{sect:doubleB}

In this section we introduce a method of constructing a virtual link diagram $\tilde D$ from a twisted link diagram $D$, called a {\it double covering diagram}, such that  the virtual link $[\tilde D]$ is the double covering of the twisted link $[D]$.  

\begin{figure}[ht]
\begin{center}
\includegraphics[width=7cm]{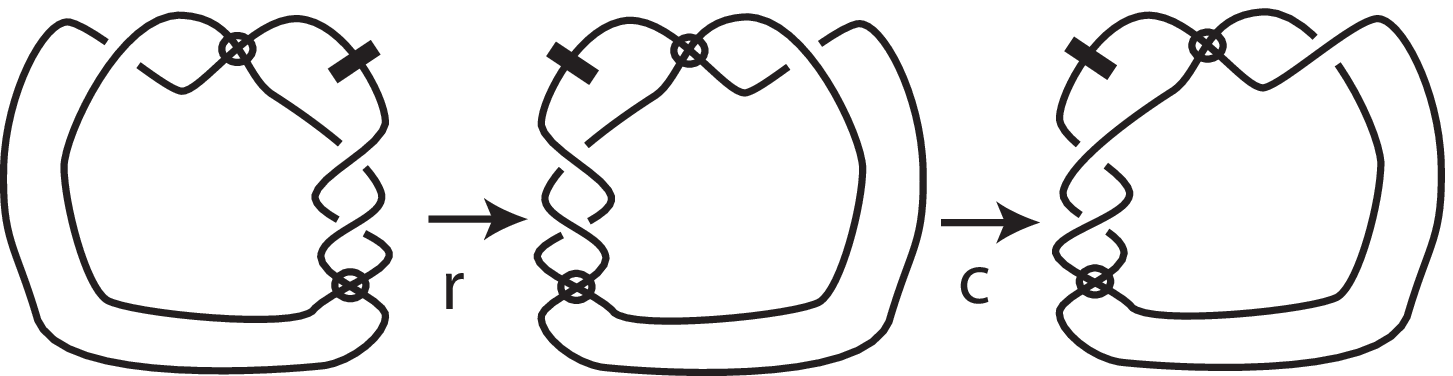}
\caption{$s  =  c \circ r : \TLD \to \TLD$}
\label{fig:figrefle} 
\end{center}
\end{figure}

Let $r : \R^2 \to \R^2, (x,y) \mapsto (-x, y)$ be the reflection along the $y$-axis and 
we denote by the same symbol $r$ the map $r : \TLD \to \TLD$ induced from $r$.  
Let 
$c : \TLD \to \TLD$ be the map changing over-under information at every real crossing, and  
let $s : \TLD \to \TLD$ be the composition $c \circ r$ of $r$ and $c$.  (Note that 
$c \circ r = r \circ c$.)  See Fig.~\ref{fig:figrefle}. 

\vspace{0.5cm}

\begin{rem}{\rm 
If $D$ is a classical link diagram, then $s(D)$ and $D$ are equivalent as classical links (and as virtual links).  On the other hand, Sawolleck \cite{rSaw} showed that there is a virtual link diagram $D$ such that $s(D)$ and $D$ are not equivalent as virtual links. 
}\end{rem}

Let $h: \R^2 \to \R, (x,y) \mapsto y$, which we call the {\it height function}. 

\begin{figure}[ht]
\begin{center}
\includegraphics[width=7cm]{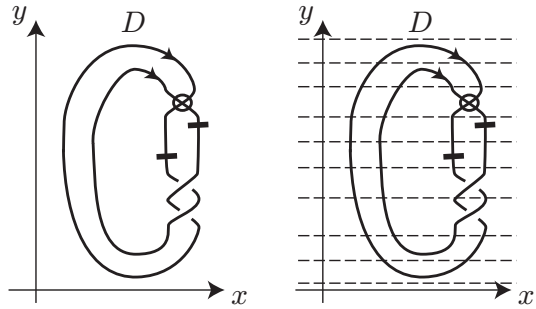}
\caption{A diagram in general position with respect to $h$}
\label{fig:figheighta} 
\end{center}
\end{figure}

We say that a twisted link diagram $D$ is {\it in general position} with respect to $h$, if (1) the critical points of $h|_D$ are minimal points or maximal points, (2) the restriction of $h$ to the set of critical points of $h|_D$, real crossings, virtual crossings, and the intersection of $D$ with bars is injective, and (3) the image of each bar is a point.  
Then there is a collection of horizontal lines of $\R^2$ such that each region 
of $\R^2$ divided by the lines contains at most one of critical points, real crossings, virtual crossings and bars. See Fig.~\ref{fig:figheighta}. 

\begin{figure}[ht]
\begin{center}
\includegraphics[width=6.5cm]{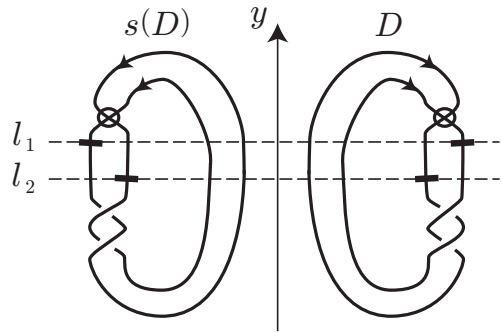}
\caption{Horizontal lines $l_1$ and $l_2$ containing the bars}
\label{fig:figheightb} 
\end{center}
\end{figure}

Let $D$ be a twisted link diagram.  
Move $D$ by an isotopy of $\R^2$ so that $D$ is in $\{ (x, y) \in \R^2 \mid x >0 \}$ and  
it is in general position with respect to $h$. 
Note that $s(D) \cup D$ is invariant under $s: \R^2 \to \R^2$.  
See Fig.~\ref{fig:figheightb}.

\begin{figure}[ht]
\begin{center}
\includegraphics[width=11cm]{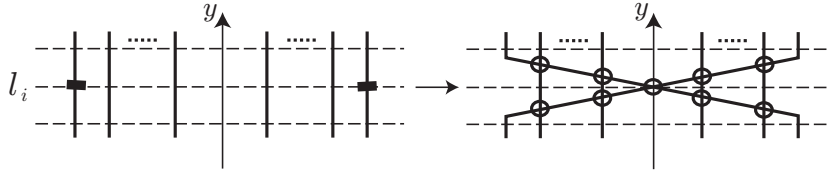}
\caption{Replacement}
\label{fig:figheightc} 
\end{center}
\end{figure}

Let $l_1, \dots,  l_n$ be horizontal lines containing the bars, where $n$ is the number of bars of $D$. 
For each $i \in \{1, \dots, n\}$, there are two bars of $s(D) \cup D$ on $l_i$, 
which are symmetric with respect to the $y$-axis. 
See Fig.~\ref{fig:figheightb}. 
Let $N(l_i)$ be a regular neighborhood of $l_i$.  There are no critical points of $h|_D$, no real crossings, no virtual crossings, and no bars in $N(l_i)$ except the two bars.  
Replace the diagram $s(D) \cup D$ in $N(l_i)$ for every $i$ as 
shown in Fig.~\ref{fig:figheightc}.  We denote by $\tilde D$ a virtual link diagram obtained this way.  

\begin{definition}{\rm 
A {\it double covering diagram} of $D$ is a diagram 
$\tilde D$ obtained above.  
}\end{definition}

\begin{figure}[ht]
\begin{center}
\includegraphics[width=5.5cm]{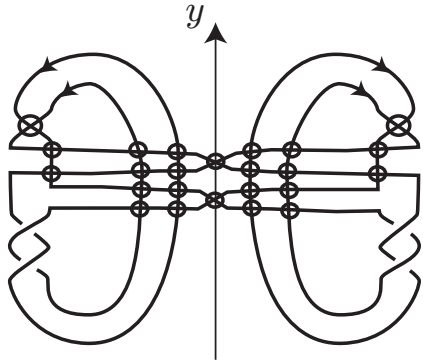}
\caption{The double covering diagram of $D$}
\label{fig:figheightd} 
\end{center}
\end{figure}

For example, when $D$ is the  diagram depicted in Fig.~\ref{fig:figheighta} then the double covering diagram is as in Fig.~\ref{fig:figheightd}. 

\begin{thm}\label{thm:doubleB}
Let $D$ be a twisted link diagram and let $\tilde D$ be a double covering diagram of $D$.  Let $K =[D] $ be the twisted link represented by $D$ and let $\tilde K = [\tilde D]$ be the virtual link represented by $\tilde D$.  Then $\tilde K$ is the double covering of $K$.  
\end{thm}

Proof.  As seen in Example~\ref{example:A}, from the abstract link diagram $A(D)$ associated with $D$, 
one can obtain an ordinary abstract link diagram  as the orientation double covering of $A(D)$ modulo inversions, 
which is also the abstract link diagram associated with $\tilde D$. \qed

\begin{thm} 
Let $D$ and $D'$ be twisted link diagrams and let 
$\tilde D$ and $\tilde D'$ be  their double covering diagrams.  
If $D$ and $D'$ are equivalent as twisted links, then 
$\tilde D$ and $\tilde D'$ are equivalent as virtual links.  
\end{thm} 

This theorem is a consequence of Theorems~\ref{thm:doubleA} and \ref{thm:doubleB}.  It is also verified by  Lemmas~\ref{lem:moveA} and \ref{lem:moveB} below.

\begin{lem}\label{lem:moveA}
Let $D$ and $D'$ be twisted link diagrams, and let 
$\tilde D$ and  $\tilde D'$ be their double covering diagrams. 
If $D$ and $D'$ are related by isotopies of $\R^2$,  V-moves and T-moves, then 
$\tilde D$ and $\tilde D'$ are related by isotopies of $\R^2$ and V-moves.  
\end{lem}

\begin{figure}[ht]
\begin{center}
\includegraphics[width=10cm]{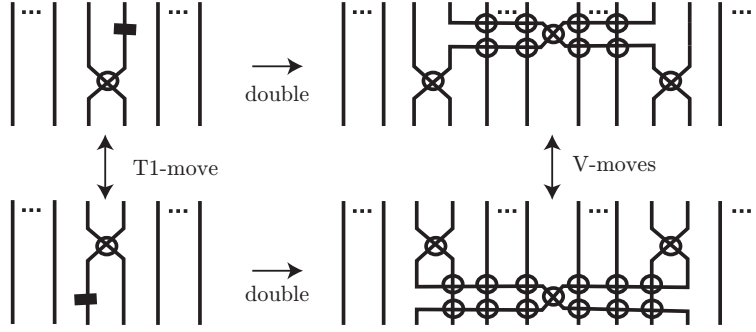}
\caption{A T1-move}
\label{fig:figheighte} 
\end{center}
\end{figure}

\begin{figure}[ht]
\begin{center}
\includegraphics[width=10cm]{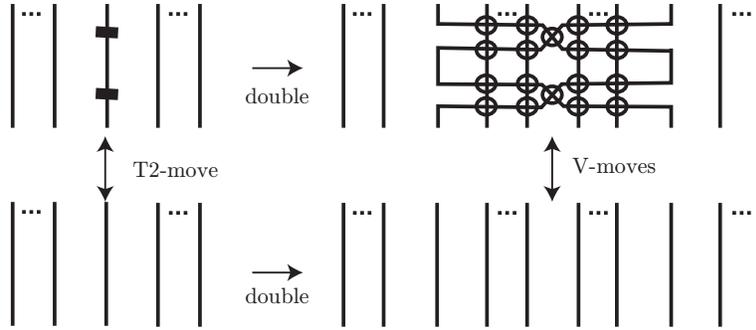}
\caption{A T2-move}
\label{fig:figheightf} 
\end{center}
\end{figure}

\begin{figure}[ht]
\begin{center}
\includegraphics[width=10cm]{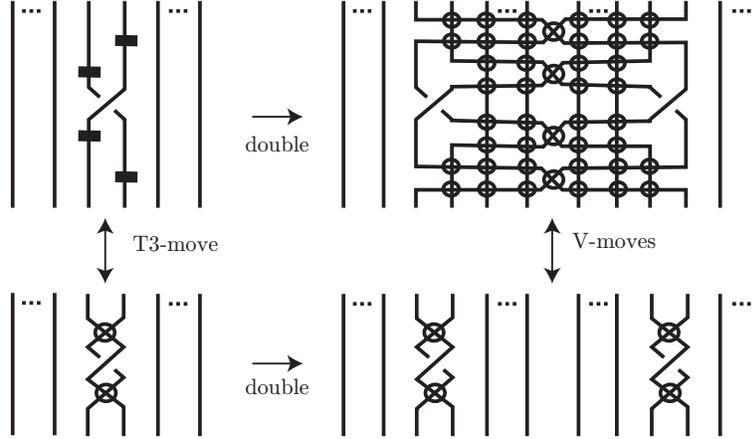}
\caption{A T3-move}
\label{fig:figheightg} 
\end{center}
\end{figure}

Proof. It is seen by considering local moves with respect to the height function $h$.  We just demonstrate three typical cases on T-moves in Fig.~\ref{fig:figheighte}, 
\ref{fig:figheightf} and \ref{fig:figheightg}.  (It is well known that two virtual link diagrams are related by isotopies of $\R^2$ and V-moves if they have the same Gauss diagram (cf. \cite{rkk, rkauD}). Using this fact, one can  easily seen that the diagrams on the right hand side of each of  these figures are related by isotopies of $\R^2$ and V-moves.)  The other cases are verified similarly and we omit them. \qed 

\vspace{0.5cm}

\begin{lem}\label{lem:moveB}
Let $D$ and $D'$ be twisted link diagrams and let 
$\tilde D$ and $\tilde D'$ be  their double covering diagrams.  
If $D$ and $D'$ are related by R-moves, then 
$\tilde D$ and $\tilde D'$ are related by R-moves.  
\end{lem}

Proof.  Assume that $D'$ is obtained from $D$ by an R-move, where we assume that $D$ is in $\{ (x, y) \in \R^2 \mid x >0 \}$ and that $D$ and $D'$ are in general position with respect to $h$.  
Let $N$ be a 2-disk where the R-move is applied to $D$.  Without of loss of generality, we may assume that $h(N)$ does not contain the images of critical points, real crossings, virtual crossings and bars of $D$ and $D'$ except the images of the crossing points involved by the R-move.  Then $\tilde D'$ is obtained from $\tilde D$ by two R-moves, one of which is applied in $N$ 
and the other is applied in $r(N)$.  \qed 

\vspace{0.5cm}
\begin{figure}[ht]
\begin{center}
\includegraphics[width=10cm]{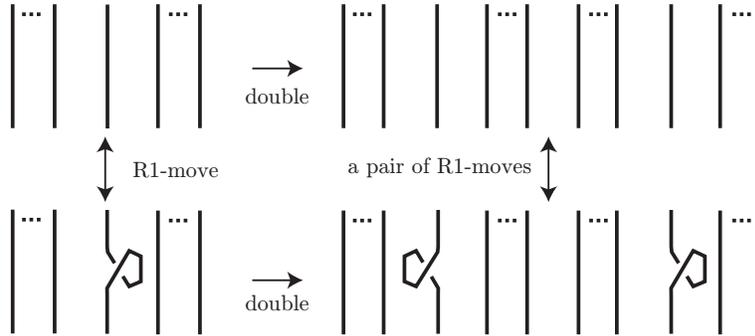}
\caption{An R1-move on $D$ changes the writhe of $\tilde D$ by $2$ or $-2$.  }
\label{fig:figheighth} 
\end{center}
\end{figure}

\begin{rem}{\rm 
As in the proof of Lemma~\ref{lem:moveB}, we see that 
an R1-move on $D$ changes the writhe of $\tilde D$ by $2$ or $-2$.  
(See Fig.~\ref{fig:figheighth}.)
}\end{rem}

\section{Links in thickened surfaces and their diagrams}\label{sect:thickenedA}

Let $F$ be a compact (not necessarily orientable) surface and let 
$I $ denote the interval $[-1,1]$. 

An {\it oriented thickened surface} over $F$ or an  {\it oriented thickening} of $F$ 
is an oriented $3$-manifold $M$ equipped with an $I$-bundle structure $F \tilde\times I$ over $F$. 
Note that the restriction of $F \tilde\times I$ 
over a connected orientable (or non-orientable) 
component of $F$ is the trivial $I$-bundle (or a twisted $I$-bundle).   We often identify $F$ with $F \times \{0\} \subset F \tilde\times I=M$.  Let $p: M=  F \tilde\times I \to F$ be the projection.  

\begin{figure}[ht]
\begin{center}
\includegraphics[width=5cm]{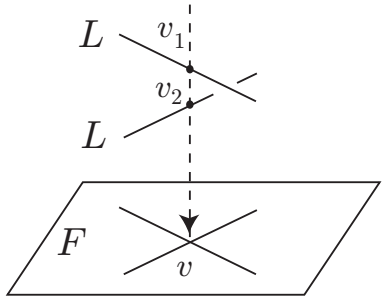}
\caption{A double point of $p(L)$ and its preimage $\{v_1, v_2\}$}
\label{fig:figdiaga} 
\end{center}
\end{figure}

For a link $L$ in an oriented thickened surface $M= F \tilde\times I$, we consider a diagram of $L$ as follows.  Move $L$ slightly so that $L$ is in general position with respect to the projection $p: M=  F \tilde\times I \to F$, i.e., $p(L)$ is an immersed circles in $F$ with (or without) transverse double points.  Let $v$ be a double point of $p(L)$ and let 
$p^{-1}(v) = \{v_1, v_2\}$.  See Fig.~\ref{fig:figdiaga}.  There is no canonical way to decide which is `over' or  `under' on $v_1$ and $v_2$. 

\vspace{0.5cm} 

\begin{figure}[ht]
\begin{center}
\includegraphics[width=7cm]{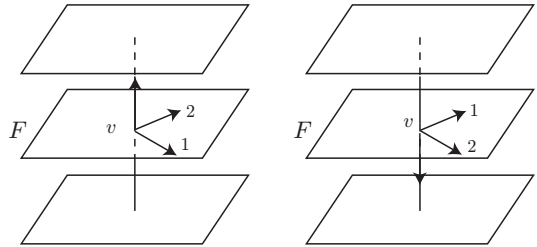}
\caption{The orientation of the fiber $I_v$ 
determined from a local orientation of $F$ at $v$}
\label{fig:figdiagb} 
\end{center}
\end{figure}

Give a local orientation of $F$ at $v$.  Using the orientation of $M= F \tilde\times I$, we can determine an orientation of the fiber $p^{-1}(v)$ as in Fig.~\ref{fig:figdiagb}. 
We call it the {\it normal orientation} of $F$ at $v$ determined by the local orientation. 

If the normal orientation points from $v_2$ toward $v_1$, then 
we say that $v_1$ is {\it over} and $v_2$ is {\it under}, and  
we call a regular neighborhood of $v_1$ (or $v_2$)  in $L$ an {\it over-crossing arc} (or an {\it under-crossing arc}) of $L$ at $v$ with respect to the normal orientation or the 
local orientation of $F$ at $v$.  The projection image in $F$ of an over-crossing arc (or an under-crossing arc) of $L$ is also called an {\it over-crossing arc} (or an {\it under-crossing arc}) of $p(L)$ at $v$.  
 (When we reverse the local orientation of $F$ at $v$, the normal orientation is reversed and the over-under information is also  reversed.)  

\vspace{0.5cm} 

Let $O$ be a collection of local orientations of $F$ at the double points of $p(L)$, which  we call a {\it local orientation system} at the double points.  
By removing a short under-arc at every double point, we obtain a {\it diagram} $D$ of $L$ in $F$.  Then we call $O$ the {\it local orientation system} for the diagram $D$.  Each double point of $p(L)$ or the corresponding point of $D$ is called a {\it crossing} of $D$. 
The image $p(L)$ is denoted by $|D|$, which is regarded as a $4$-valent graph in $F$ and called the {\it underlying  $4$-valent graph} of $D$.  
We call the pair $(D, O)$ or the triple $(F, D, O)$ 
a {\it generalized link diagram} of $L$.  

\begin{figure}[ht]
\begin{center}
\includegraphics[width=7cm]{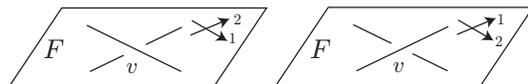}
\caption{An inversion at $v$}
\label{fig:figdiagc} 
\end{center}
\end{figure}

Note that the diagram $D$ depends on a choice of $O$. If we reverse the local orientation at a crossing $v$, then the over-under information at $v$ is reversed.  See Fig.~\ref{fig:figdiagc}.  

\vspace{0.5cm} 

Two generalized link diagrams $(D, O)$ and $(D', O')$ on $F$ are 
said to be {\it isotopic} if there is an isotopy of $F$ carrying $|D|$ to $|D'|$ such that for each crossing point of $D$, the local orientation and the over-under information are preserved.  (We usually do not distinguish isotopic  generalized link diagrams.) 

Let $(D, O)$ be a generalized link diagram on $F$ 
and let $v$ be a crossing point of $D$. 
A generalized link diagram $(D', O')$ is said to be obtained from $(D, O)$ 
by an {\it inversion} at $v$ if it is obtained by reversing the local orientation and 
the over-under information at $v$. See Fig.~\ref{fig:figdiagc}. 

A generalized link diagram $(D', O')$ is said to be obtained from $(D, O)$ 
by a {\it Reidemeister move} if there exists an oriented $2$-disk $N$ in $F$ such that 
the diagram $D'$ is obtained from $D$ by a Reidemeister move  in 
$N$ and the local orientations of the crossings of $D$ and $D'$ in $N$ are compatible with the orientation of $N$.  

\begin{definition}{\rm 
Two generalized link diagrams $(D, O)$ and $(D', O')$ on $F$ are 
 {\it R-equivalent} on $F$  if there exists a finite sequence of generalized link diagrams on $F$, 
$(D, O)= (D_0, O_0),$  $(D_1, O_1), \dots,$  $(D_n, O_n) = (D', O')'$ such that for each $i \in \{1, \dots, n\}$, $(D_i, O_i)$ is obtained from $(D_{i-1}, O_{i-1})$ by an isotopy of $F$, 
an inversion or a Reidemeister move.   
}\end{definition}

\begin{thm}\label{thm:LTSisotopic}
Let $L$ and $L'$ be links in a thickened surface $M= F \tilde\times I$, and let 
$(D, O)$ and $(D', O')$  be their generalized link diagrams.  The links $L$ and $L'$ are isotopic in $M$ if and only if $(D, O)$ and $(D', O')$ are R-equivalent on $F$. 
\end{thm}

Proof.  The if part is obvious.  Suppose that $L$ and $L'$ are isotopic.  Let $\cal U$ be an open covering of $F$ such that each element $U \in \cal U$ is an open $2$-disk in $F$ and the $I$-bundle structure of $M$ over $U$ is trivial.  Without loss of generality we may assume that $L$ and $L'$ are PL  links in $M$ and there is a finite sequence of links $L=L_0, L_1, \dots, L_n=L'$ such that for each $i$, $L_i$ is in general position with respect to the projection $M \to F$ and $L_i$ is obtained from $L_{i-1}$ by a move along a $2$-simplex $\Delta_i$ in $p^{-1}(U)$ for some $U \in \cal U$.  Let $(D_0, O_0), (D_1, O_1), \dots, (D_n, O_n)$ be the sequence of generalized link diagrams such that 
$(D, O) = (D_0, O_0)$, and for each $i \in \{1, \dots, n\}$, $(D_i, O_i)$ is a generalized link diagram of $L_i$ such that it is obtained from $(D_{i-1}, O_{i-1})$ by a finite sequence of Reidemeister moves and isotopies of $F$.  Since $(D_n, O_n)$ and $(D', O')$ are generalized link diagrams of the same link $L$, they are related by inversions. \qed

In what follows, a pair $(M, L)$ of an oriented thickened surface $M$ and a link $L$ in $M$ is also referred to as a link in an oriented thickened surface.   For a generalized link diagram $(D, O)$ on a surface $F$, the triple $(F, D, O)$ is also referred to as a generalized link diagram.  

\begin{definition}{\rm 
Two links in oriented thickened surfaces $(M, L)$ and $(M', L')$ are 
{\it equivalent} if there is an orientation-preserving and fiber-preserving homeomorphism from $M$ to $M'$ that sends the link $L$ onto a link isotopic $L'$ in $M'$.  
}\end{definition} 

\begin{definition}{\rm 
Two generalized link diagrams $(F, D, O)$ and $(F', D', O')$ are {\it homeomorphic} if there is a homeomorphism from $F$ to $F'$ that sends $(D, O)$ to $(D', O')$.  
}\end{definition} 

\begin{definition}{\rm 
Two generalized link diagrams $(F, D, O)$ and $(F', D', O')$ are {\it equivalent} if there is a 
finite sequence of generalized link diagrams $(F, D, O)= (F_0, D_0, O_0), (F_1, D_1, O_1),$  $\dots,$ 
$(F_n, D_n, O_n)=(F', D', O')$ such that for each $i$, (1) $(F_i, D_i, O_i)$ 
is homeomorphic to $(F_{i-1}, D_{i-1}, O_{i-1})$ or (2) $F_i = F_{i-1}$ and 
$(D_i, O_i)$ is R-equivalent to $(D_{i-1}, O_{i-1})$ on $F_i = F_{i-1}$. 
}\end{definition} 

\begin{thm}\label{thm:GLDequiv}
Let $(M, L)$ and $(M', L')$ be links in oriented thickened surfaces 
over $F$ and $F'$, respectively. 
Let $(F, D, O)$ and $(F', D', O')$ be generalized link diagrams of 
$(M, L)$ and $(M', L')$.  
The links $(M, L)$ and $(M', L')$ are equivalent if and only if 
$(F, D, O)$ and $(F', D', O')$  are equivalent. 
\end{thm}

Proof. 
Suppose that $(M, L)$ and $(M', L')$ are equivalent.  By definition, there is an orientation-preserving and fiber-preserving homeomorphism from $M$ to $M'$ that sends $L$ onto a link $L''$ which is isotopic to $L'$ in $M'$.  Then $L''$ has a generalized link diagram $(F', D'', O'')$ that is homeomorphic to $(F, D, O)$.  On the other hand, by Theorem~\ref{thm:LTSisotopic}, 
$(F', D', O')$ is R-equivalent to $(F', D'', O'')$.  Thus we have the only if part.  The if  part is obvious.  \qed

\section{Stable equivalence on links in thickened surfaces}\label{sect:thickenedB}

We introduce an equivalence relation on links in oriented thickened surfaces, 
called  the stable equivalence relation.  Our definition of the stable equivalence relation is stated in a different way from those given in \cite{rBourg} and \cite{rCKS}, although they are the same equivalence relation.  Using our definition, one can easily understand the relationship between 
abstract links (or twisted links) and stable equivalence classes of  links in oriented  thickened surfaces.

Let $M$ be an oriented thickened surface over $F$ and let $p : M = I \tilde\times F \to F$ be the projection.  Let $c$ be a  $2$-sided  simple loop in $F$.  A regular neighborhood $A$ of $c$ is an annulus in $F$. Remove ${\rm int} A$, the interior of $A$,  from $F$ and attach two $2$-disks along the boundary. We obtain a surface $F'$, which is called the surface obtained from $F$ by a $2$-handle surgery along $c$. An oriented thickened surface over $F'$ is obtained by removing $p^{-1}({\rm int} A)$ from $M$ and attaching two trivial $I$-bundles over $2$-disks along $p^{-1}(\partial A)$.  We call it the oriented thickened surface obtained from $M$ by a {\it $2$-handle surgery}  along $c$.   The inverse operation is called a {\it $1$-handle surgery}. 

Let $(M, L)$ be a link in an oriented thickened surface, and let 
$p : M= I \tilde\times F \to F$  be the projection.  
Let $c$ be a $2$-sided simple loop in $F$ with a regular neighborhood $A$ in $F$ which is disjoint from $p(L)$. Let $M'$ be the oriented thickened surface obtained from $M$ by a $2$-handle surgery along $c$, and let $L'$ be the link $L$ in $M'$. Then we say that 
$(M', L')$ is obtained from $(M,L)$ by a {\it $2$-handle ambient surgery}, and that 
$(M, L)$ is obtained from $(M',L')$ by a {\it $1$-handle ambient surgery}.  

Let $(M, L)$ be a link in an oriented thickened surface, and let 
$p : M= I \tilde\times F \to F$  be the projection.  Suppose that $F$ has a connected component, say $F_0$, such that $p(L)$ is disjoint from $F_0$.  We call the connected component of $M$ over $F_0$ a {\it nugatory ambient component}.  
Let $M'$ be the oriented thickened surface obtained from $M$ by removing a  nugatory ambient component, and let $L'$ be the link $L$ in $M'$.  Then we say that 
$(M', L')$ is obtained from $(M,L)$ by  {\it elimination of a nugatory ambient component}, and that 
$(M, L)$ is obtained from $(M',L')$ by   {\it addition of a nugatory ambient  component}.

\begin{definition}\label{def:stableA}{\rm 
Two links in oriented thickened surfaces $(M, L)$ and $(M, L')$ are 
{\it stably equivalent} if they are related by a finite sequence of links in oriented thickened surfaces, 
$(M, L)= (M_0, L_1), (M_1, L_1), \dots, (M_n, L_n)=(M', L')$ such that 
for each $i \in \{1, \dots, n\}$, $(M_i,L_i)$ is equivalent to $(M_{i-1}, L_{i-1})$ or 
$(M_i,L_i)$ is obtained from $(M_{i-1}, L_{i-1})$ by a $1$-handle ambient surgery, a $2$-handle ambient surgery, elimination of a nugatory ambient component, or addition of a nugatory ambient component.
}\end{definition}

Let $(F, D, O)$ be a generalized link diagram. 
Let $c$ be a $2$-sided simple loop in $F$ with a regular neighborhood $A$ in $F$ which is disjoint from $D$.  Let $F'$ be the surface obtained from a $2$-handle surgery along $c$ with $A$ and let $D'$ and $O'$ be $D$ and $O'$ in $F'$.  Then we say that $(F', D', O')$ is obtained from $(F, D, O)$ by 
a {\it $2$-handle ambient surgery}, and that 
$(F, D, O)$ is obtained from $(F', D', O')$ by 
a {\it $1$-handle ambient surgery}.  

Let $(F, D, O)$ be a generalized link diagram.  Suppose that $F$ has a connected component, say $F_0$, such that $D$ is disjoint from $F_0$.  We call  $F_0$ a {\it nugatory ambient component}.  
Let $F'$ be the  surface obtained from $F$ by removing a nugatory ambient component, and let $D'$ and $O'$ be $D$ and $O'$ in $F'$. 
Then we say that 
$(F', D', O')$ is obtained from $(F, D, O)$ by  {\it elimination of a nugatory ambient component}, and that 
$(F, D, O)$ is obtained from $(F', D', O')$ by   {\it addition of a nugatory ambient  component}.  

\begin{definition}\label{def:stableB}{\rm 
Two generalized link diagrams $(F, D, O)$ and $(F', D', O')$ are 
{\it stably equivalent} if they are related by a finite sequence of 
generalized link diagrams, 
$(F, D, O)= (F_0, D_0, O_0), (F_1, D_1, O_1), \dots, (F_n, D_n, O_n)=(F', D', O')$ such that 
for each $i \in \{1, \dots, n\}$, $(F_i, D_i, O_i)$ is equivalent to  $(F_{i-1}, D_{i-1}, O_{i-1})$ or 
$(F_i, D_i, O_i)$ is obtained from $(F_{i-1}, D_{i-1}, O_{i-1})$  by a $1$-handle ambient surgery, a $2$-handle ambient surgery, elimination of a nugatory ambient component, or addition of a nugatory ambient component.
}\end{definition}

\begin{rem}{\rm 
Our definitions of stable equivalences (Definitions~\ref{def:stableA} and \ref{def:stableB}) involve stabilizations producing  (a thickening of) a non-orientable surface, and destabilizations producing (a thickening of) an orientable surface.  For example, let $(F, D, O)$ be a generalized link diagram such that $F$ is a projective plane. Let $c_0$ be a $1$-sided simple loop in $F$ and suppose that there is a regular neighborhood $N$ of $c_0$ which is disjoint from $D$. Let $c$ be a $2$-sided simple loop in $F$ which is in $N$ and is parallel to $\partial N$. Apply a $2$-handle ambient surgery along $c$ to $(F, D, O)$ and we obtain a generalized link diagram, say $(F', D', O')$, such that  $F'= F'_1 \cup F'_2$ where $F'_1$ is a $2$-sphere and $F'_2$ is a projective plane and  that $D'$ is disjoint from $F'_2$.  We can remove $F'_2$ from the ambient surface by elimination of a nugatory ambient component. Then we have a diagram in the $2$-sphere $F'_1$. 
}\end{rem}

\begin{thm}\label{thm:GLDstableequiv}
Let $(M, L)$ and $(M, L')$ be links in oriented thickened surfaces over $F$ and $F'$, respectively.  Let 
$(F, D, O)$ and $(F', D', O')$ be generalized link diagrams of $(M, L)$ and $(M, L')$.  
The links $(M, L)$ and $(M, L')$ are stably equivalent if and only if $(F, D, O)$ and $(F', D', O')$ are stably equivalent. 
\end{thm}

Proof.  It follows from Theorem~\ref{thm:GLDequiv} and the definitions of 
a $1$-handle ambient surgery, a $2$-handle ambient surgery, elimination of a nugatory ambient component, and addition of a nugatory ambient component. 
\qed 

\begin{thm}[cf. \cite{rBourg}]
There is a bijection between the set of stable equivalence classes of generalized link diagrams over closed surfaces (or equivalently the set of stable equivalence classes of links in oriented thickened surfaces over closed surfaces by Theorem~\ref{thm:GLDstableequiv}) and the set of abstract links. 
\end{thm} 

This theorem is essentially due to Bourgoin \cite{rBourg}, although our definition of the stable equivalence relation is stated differently.   Refer also to  
\cite{rCKS, rFRS, rkk} for the relationships among virtual links, (ordinary) abstract links, and stable equivalence classes of link diagrams on closed oriented surfaces. 

\vspace{0.5cm} 

Proof. 
Let $(F, D, O)$ be a generalized link diagram with $F$ closed.  Let $N(D)$ denote a regular neighborhood of the underlying $4$-valent graph $|D|$ of $D$ in $F$. Then we have an abstract link diagram $(N(D), D, O)$.  It is easily seen that if $(F, D, O)$ and $(F', D', O')$ are stably equivalent, then $(N(D), D, O)$ and $(N(D'), D', O')$ are equivalent as abstract links.  Thus we have a well-defined map 
$$ J : 
\left\{ 
\begin{array}{c}
\mbox{ stable equivalence classes of generalized }  \\
\mbox{ link diagrams over closed surfaces }
\end{array}
\right\} 
\to \{ \mbox{ abstract links } \}$$ 
by sending the class of  $(F, D, O)$ to the class of $(N(D), D, O)$.  

Let $(\Sigma, D, O)$ be an abstract link diagram.  Let $E$ be a compact surface such that $\partial \Sigma$ and $\partial E$ has the same number of connected components.  Take a homeomorphism $f : \partial E \to \partial \Sigma$ and consider a closed surface $F$ that is the union of $E$ and $\Sigma$ by identifying their boundaries with $f$.  Then we have a generalized link diagram $(F, D, O)$ with $F$ closed. We call it a {\it closure} of $(\Sigma, D, O)$.  
Let $(F', D, O)$ be another closure of $(\Sigma, D, O)$ obtained by using a compact surface $E'$ and a homeomorphism $f' : \partial E' \to \partial \Sigma$.  
Then $(F, D, O)$ and $(F', D, O)$ are stably equivalent. (For  
each of $(F, D, O)$ and 
$(F', D, O)$, apply $2$-handle ambient surgeries along the loops $\partial \Sigma$ followed by elimination of all nugatory ambient components.  Then the results are  the same.)  One can easily verified that if $(\Sigma, D, O)$ and $(\Sigma', D', O')$ are equivalent as abstract links, then $(F, D, O)$ and $(F', D', O')$ are stably equivalent. 
Thus   we have a well-defined map 
$$ K : 
\{ \mbox{ abstract links } \} \to 
\left\{ 
\begin{array}{c}
\mbox{ stable equivalence classes of generalized }  \\
\mbox{ link diagrams over closed surfaces }
\end{array}
\right\} 
$$ 
by sending the class of  $(\Sigma, D, O)$ to the class of a closure of  $(\Sigma, D, O)$.  Since $J$ and $K$ are the inverse maps each other, we have the result. \qed 

Kuperberg \cite{rkup} proved that for a stable equivalence class of links in oriented  thickened surfaces over oriented closed surfaces there is a unique representative with the surface having minimal genus.  Bourgoin \cite{rBourg} generalized it for  a stable equivalence class of links in oriented thickened surfaces over closed surfaces.

\section{Twisted link groups and quandles}\label{sect:groups}

Bourgoin \cite{rBourg}   introduced the twisted link group 
$\tilde \Pi(D)$ for a twisted link diagram $D$, and proved that it is an invariant of a twisted link.  
He proposed a problem asking a geometric interpretation of the twisted link group.   We answer it in this section.  We also give a geometric interpretation of the twisted link quandle defined in \cite{rkamF}.

Let $D$ be a twisted link diagram, and 
let $a_1, \dots, a_n$ be the  semi-arcs of $D$.  Here a {\it semi-arc} of $D$ means an arc of the simple or non-simple arcs obtained by cutting $D$ at all real crossings and at all bars.  (For simplicity, we number the semi-arcs $a_1, \dots, a_n$  along the orientation.)  
For each semi-arc $a_i$ of $D$, assign two labels $x_i$ and $y_i$. 
These are generators of the group $\tilde \Pi(D)$.  

\vspace{0.5cm} 
\begin{figure}[ht]
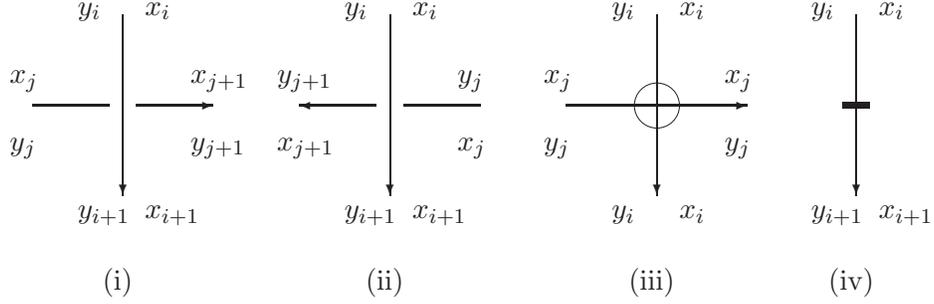

\begin{center}
\begin{tabular}{cccc}
\fgknotgpp{.6mm}&\fgknotgpm{.6mm}&\fgknotgpv{.6mm}&\fgknotgpt{.6mm}\\
 & & & \\ 
(i)  & (ii)  & (iii) & (iv) \\
\end{tabular}
\caption{Labels on arcs around real crossings, virtual crossings and bars}\label{fig:labels}
\end{center}
\end{figure}

For a positive crossing $c$ with labels as in Fig.~\ref{fig:labels} (i), we consider 4 relations:   
$$ 
x_{i+1} = x_i \quad 
x_{j+1} = x_i^{-1} x_j x_i, \quad 
y_{i+1} = y_j^{-1} y_i y_j \quad 
\mbox{and} \quad 
y_{j+1}= y_j. $$ 
The first two relations $\{ x_{i+1} = x_i,  x_{j+1} = x_i^{-1} x_j x_i \}$ is denoted by $\mbox{upper-rel}(c)$ and the other two relations $\{ y_{i+1} = y_j^{-1} y_i y_j, y_{j+1}= y_j \}$ is denoted by $\mbox{lower-rel}(c)$.  

For a negative crossing $c$ with labels as in Fig.~\ref{fig:labels} (ii), we consider 4 relations:   
$$ 
x_{i+1} = x_i \quad 
x_{j+1} = x_i x_j x_i^{-1}, \quad 
y_{i+1} = y_j y_i y_j^{-1} \quad 
\mbox{and} \quad 
y_{j+1}= y_j. $$ 
The first two relations $\{ x_{i+1} = x_i,  x_{j+1} = x_i x_j x_i^{-1} \}$ is denoted by $\mbox{upper-rel}(c)$ and the other two relations $\{ y_{i+1} = y_j y_i y_j^{-1}, y_{j+1}= y_j \}$ is denoted by $\mbox{lower-rel}(c)$.  

For a bar $b$ with labels as in Fig.~\ref{fig:labels} (iv), we consider 2 relations:  
$$ x_{i+1} = y_i \quad \mbox{and} \quad 
y_{i+1} = x_i. $$
We denote these two relations by $\mbox{bar-rel}(b)$.  

The {\it twisted link group} of $D$,  denoted by $\tilde \Pi(D)$,  is the group determined by a group presentation whose generators are the symbols  $x_i$ and $y_i$ $(i=1, \dots, n)$ and the defining relations are the 4 relations for each crossing and the 2 relations for each bar.  

Unless otherwise stated, the virtual link group means the upper virtual link group (cf. \cite{rkk, rkauD}). 

\begin{thm}\label{thm:twistedgroupA}
The twisted link group of $D$ is the virtual link group of $\tilde D$.  
\end{thm}  

Proof. Without loss of generality we may assume that $D$ is in $\{ (x, y) \in \R^2 \mid x >0 \}$ and in general position with respect to the height function $h$.  Let $a_1, \dots, a_n$ be the  semi-arcs of $D$.  
By definition, the twisted link group $\tilde \Pi(D)$ of $D$ has a presentation  whose generators are 
$x_1, \dots, x_n, y_1, \dots,  y_n$ and the defining relations are 
$$ 
\mbox{upper-rel}(c), ~  \mbox{lower-rel}(c) \quad (c : \mbox{crossings of $D$}), \quad 
\mbox{and}
\quad \mbox{bar-rel}(b) \quad (b : \mbox{bars of $D$}). $$

Let $D^\vee$ and $s(D)^\vee$ be the 
diagrams obtained from $D$ and $s(D)$ by cutting arcs at all bars.  
(Although the diagram $D^\vee$ has some end-points,  the (upper) virtual link group can be defined even for such a diagram.)  
The virtual link group $\Pi(D^\vee)$ 
of $D^\vee$ has a presentation 
$$\langle x_1, \dots, x_n \mid \mbox{upper-rel}(c) \quad (c : \mbox{crossings of $D$}) \rangle$$
and the virtual link group $\Pi(s(D)^\vee)$ 
of $s(D)^\vee$ has a presentation 
$$\langle y_1, \dots,  y_n \mid \mbox{lower-rel}(c) \quad (c : \mbox{crossings of $D$}) \rangle.$$ 
The diagram $\tilde D$ is obtained from the union of $D^\vee$ and $s(D)^\vee$ by 
connecting their endpoints with some arcs, and they make additional relations  
$\mbox{bar-rel}(b)$ for all bars $b$.  
Thus the virtual link group $\Pi(\tilde D)$ is the amalgamated free product 
of $\Pi(D^\vee)$ and $\Pi(s(D)^\vee)$ with amalgamation by using $\mbox{bar-rel}(b)$ for all  bars $b$.  This coincides with the twisted link group $\tilde \Pi(D)$. \qed 

Now we give a geometric interpretation of the twisted link group.  

\begin{thm}\label{thm:twistedgroupB}
Let $D$ be a twisted link diagram and let $\tilde D$ be a double covering diagram of $D$. 
Let $\tilde L$ be a link in the oriented thickened surface $\tilde\Sigma\times I$ presented by the ordinary 
abstract link diagram $A(\tilde D) =(\tilde\Sigma, A(\tilde D))$.  Then 
the twisted link group $\tilde\Pi (D)$ is isomorphic to the fundamental group 
of the complement of $\tilde L$ in the singular $3$-manifold $\tilde\Sigma\times I / \tilde\Sigma\times\{1\}$.  
\end{thm}

Proof. 
By Theorem~\ref{thm:twistedgroupA}, the twisted link group $\tilde\Pi (D)$ is isomorphic to the virtual link group $\Pi(\tilde D)$ of $\tilde D$.  Since $\tilde D$ is a virtual link diagram, 
it is shown in \cite{rkk} that the virtual link group $\Pi(\tilde D)$  is isomorphic to 
the fundamental group of the complement of $\tilde L$ in $\tilde\Sigma\times I  / \tilde\Sigma\times\{1\}$. 
\qed 

\begin{prop}\label{thm:twistedgroupC}
Let  $\tilde D$ be a double covering diagram of a twisted link diagram $D$. 
The  (upper) virtual link group $\Pi(\tilde D)$ is isomorphic to the lower 
virtual link group $\Pi^{l}(\tilde D)$.  
\end{prop} 

Proof. 
We may assume that $D$ is in $\{ (x, y) \in \R^2 \mid x >0 \}$ and in general position with respect to the height function $h$.  Let $a_1, \dots, a_n$ be the  semi-arcs of $D$.  
Let $D^\vee$ and $s(D)^\vee$ be the 
diagrams obtained from $D$ and $s(D)$ by cutting arcs at all bars.  

The lower virtual link group $\Pi^{l}(D^\vee)$ 
of $D^\vee$ has a presentation 
$$\langle y_1, \dots,  y_n \mid \mbox{lower-rel}(c) \quad (c : \mbox{crossings of $D$}) \rangle$$ 
and the lower virtual link group $\Pi^{l}(s(D)^\vee)$ 
of $s(D)^\vee$ has a presentation 
$$\langle x_1, \dots, x_n \mid \mbox{upper-rel}(c) \quad (c : \mbox{crossings of $D$}) \rangle.$$
The lower virtual link group $\Pi^{l}(\tilde D)$ is the amalgamated free product 
of $\Pi^{l}(D^\vee)$ and $\Pi^{l}(s(D)^\vee)$ with amalgamation by using $\mbox{bar-rel}(b)$ for all  bars $b$.  
Thus $\Pi^{l}(\tilde D)$ has the same group presentation with $\Pi(\tilde D)$. \qed 

\begin{cor}
A virtual link $K$ is not the double covering of a twisted link  if  the virtual link group of $K$ is not isomorphic to the lower virtual link group of $K$.  
\end{cor} 

The first author introduced the twisted link quandle of a twisted link (\cite{rkamF}). For a twisted link diagram $D$, the {\it twisted link quandle} $\tilde Q(D)$ is the quandle determined by a quandle presentation obtained from the group presentation of the definition of $\tilde\Pi(D)$ by changing each conjugate $a^{-1}ba$ by the quandle operation $b \ast a$ and $aba^{-1}$ by the dual operation $b \overline\ast a$. 
Then the twisted link quandle is an invariant of a twisted link  (\cite{rkamF}).  

\begin{thm}\label{thm:twistedquandleA}
The twisted link quandle of $D$ is the (upper) virtual link quandle of $\tilde D$.  
\end{thm}  

Proof. 
By the same argument with the proof of Theorem~\ref{thm:twistedgroupA}, we see the result.  
\qed 

\begin{thm}\label{thm:twistedquandleB}
Let $D$ be a twisted link diagram and let $\tilde D$ be a double covering diagram of $D$. 
Let $\tilde L$ be a link in the oriented thickened surface $\tilde\Sigma\times I$ presented by the ordinary 
abstract link diagram $A(\tilde D) =(\tilde\Sigma, A(\tilde D))$.  Then 
the twisted link quandle $\tilde Q(D)$ is isomorphic to the fundamental  quandle 
of $\tilde L$ in the singular $3$-manifold  $\tilde\Sigma\times I / \tilde\Sigma\times\{1\}$.  
\end{thm}

Proof. 
By Theorem~\ref{thm:twistedquandleA}, 
 the twisted link quandle $\tilde Q(D)$ is isomorphic to the (upper) virtual link quandle $Q(\tilde D)$ of $\tilde D$.  Since $\tilde D$ is a virtual link diagram, 
it is shown in \cite{rkk} that the virtual link quandle $Q(\tilde D)$  is isomorphic to 
the fundamental quandle of $\tilde L$ in 
$\tilde\Sigma\times I  / \tilde\Sigma\times\{1\}$. 
\qed 

By the same argument with the proof of Proposition~\ref{thm:twistedgroupC}, we see the following. 

\begin{prop}\label{thm:twistedquandleC}
Let  $\tilde D$ be a double covering diagram of a twisted link diagram $D$. 
The  (upper) virtual link quandle $Q(\tilde D)$ is isomorphic to the lower 
virtual link quandle $Q^{l}(\tilde D)$.  
\end{prop} 

\begin{cor}
A virtual link $K$ is not the double covering of a twisted link  if  the virtual link quandle of $K$ is not isomorphic to the lower virtual link quandle of $K$.  
\end{cor} 

Bourgoin \cite{rBourg} proved that if $D$ is a virtual link diagram, then 
$\tilde \Pi (D)$ is the free product of the (upper) virtual link group $\Pi (D)$ and the lower virtual link group $\Pi^{l}(D)$.  This fact is also seen from the proof of Theorem~\ref{thm:twistedgroupA}.  Similarly, from the proof of Theorem~\ref{thm:twistedquandleA}, we obtain the following. 

\begin{prop}
If $D$ is a virtual link diagram, then 
$\tilde Q(D)$ is the free product of the (upper) virtual link quandle  $Q(D)$ and the lower virtual link quandle $Q^{l}(D)$. 
\end{prop}

 \end{document}